\documentclass{IEEEtran4PSCC}
\ifCLASSINFOpdf

\else

\fi
\usepackage[cmex10]{amsmath}

\usepackage{amssymb}

\usepackage{graphics,graphicx,epic,eepic,epsfig,times,units} 

\usepackage{float}

\usepackage{psfrag}

\newcommand{\reals}{{\mbox{\bf R}}}
\newcommand{\SOC}{{\mbox{\bf L}}}

\usepackage{ifpdf}

\usepackage{cite}

\usepackage{siunitx}
\DeclareSIUnit \VAr {VAr} 
\DeclareSIUnit \VA {VA} 
\DeclareSIUnit \rad {Radians} 

\usepackage{multirow}

\usepackage{cases}

\usepackage{empheq}

\usepackage{cleveref}
\crefformat{equation}{(#2#1#3)}
\crefrangeformat{equation}{(#3#1#4)--(#5#2#6)}
\crefmultiformat{equation}{(#2#1#3)}%
{,~(#2#1#3)}{, (#2#1#3)}{,~(#2#1#3)}

\usepackage{eqnarray}

\usepackage{xcolor}

\usepackage{tabularx}

\usepackage[ruled,vlined]{algorithm2e}

\newenvironment{Model}[1][htb]
  {
   \begin{algorithm}[#1]%
  }{\end{algorithm}}

\begin{document}

\title{Tight LP Approximations for the Optimal Power Flow Problem}

\author{
\IEEEauthorblockN{Sleiman~Mhanna\\Gregor~Verbi\v{c}\\Archie~C.~Chapman}
\IEEEauthorblockA{School of Electrical and Information Engineering, \\
The University of Sydney\\
Sydney, Australia\\
\{sleiman.mhanna,gregor.verbic,archie.chapman\}@sydney.edu.au}
}

\maketitle

\begin{abstract}
DC power flow approximations are ubiquitous in the electricity industry. However, these linear approximations fail to capture important physical aspects of power flow, such as the reactive power and voltage magnitude, which are crucial in many applications to ensure voltage stability and AC solution feasibility. This paper proposes two LP approximations of the AC optimal power flow problem, founded on tight polyhedral approximations of the SOC constraints, in the aim of retaining the good lower bounds of the SOCP relaxation and relishing the computational efficiency of LP solvers. The high accuracy of the two LP approximations is corroborated by rigorous computational evaluations on systems with up to $9241$ buses and different operating conditions. The computational efficiency of the two proposed LP models is shown to be comparable to, if not better than, that of the SOCP models in most instances. This performance is ideal for MILP extensions of these LP models since MILP is computationally more efficient than MIQCP.
\end{abstract}

\section*{Notation}
\addcontentsline{toc}{section}{Notation}
\subsection{Input data and operators}
\begin{IEEEdescription}[\IEEEsetlabelwidth{$P_{i}^{g}/Q_{i}^{g}$}\IEEEusemathlabelsep]
\item[$\mathcal{B}$] Set of buses in the power network.
\item[$\mathcal{B}_{i}$] Set of buses connected to bus $i$.
\item[$b^\text{sh}_{ij}$] Shunt susceptance (p.u.) in the $\pi$-model of line $ij$.
\item[$c0^{g}_{i}$] Constant coefficient ($\SI{}{\$}$) term of generator $g$'s cost function.
\item[$c1^{g}_{i}$] Coefficient ($\SI{}{\$\per\mega\watt}$) of the linear term of generator $g$'s cost function.
\item[$c2^{g}_{i}$] Coefficient ($\SI{}{\$\per\mega\watt\squared}$) of the quadratic term of generator $g$'s cost function.
\item[$\mathcal{G}$] Set of all generators $(g,i)$ in the power network such that $g$ is the generator and $i$ is the bus connected to it.
\item[$\mathrm{i}$] Imaginary unit.
\item[$\mathcal{L}$] Set of all transmission lines $ij$ where $i$ is the ``from'' bus.
\item[$\mathcal{L}_{t}$] Set of all transmission lines $ij$ where $i$ is the ``to'' bus.
\item[$P_{i}^{d}/Q_{i}^{d}$] Active/reactive power demand ($\SI{}{\mega\watt}/\SI{}{\mega\VAr}$) at bus $i$.
\item[$\overline{S}_{ij}$] Apparent power rating ($\SI{}{\mega\VA}$) of line $ij$.
\item[$\underline{\theta}_{ij}^{\Delta}$] Lower limit of the difference of voltage angles at buses $i$ and $j$.
\item[$\overline{\theta}_{ij}^{\Delta}$] Upper limit of the difference of voltage angles at buses $i$ and $j$.
\item[$\theta_{ij}^{\text{shift}}$] Phase shift ($\SI{}{\rad}$) of phase shifting transformer connected between buses $i$ and $j$ ($\theta_{ij}^{\text{shift}}=0$ for a transmission line).
\item[$\tau_{ij}$] Tap ratio magnitude of phase shifting transformer connected between buses $i$ and $j$ ($\tau_{ij}=1$ for a transmission line).
\item[$y_{ij}$] Series admittance (p.u.) in the $\pi$-model of line $ij$.
\item[$\Im\left\{\bullet\right\}$] Imaginary value operator.
\item[$\Re\left\{\bullet\right\}$] Real value operator.
\item[$\underline{\bullet}/\overline{\bullet}$] Minimum/maximum magnitude operator.
\item[$\left|\bullet\right|$] Magnitude operator/Cardinality of a set.
\item[$\bullet^*$] Conjugate operator.
\item[$\bullet \succeq $] Matrix inequality sign in the positive semidefinite sense.
\end{IEEEdescription}

\subsection{Decision variables}
\begin{IEEEdescription}[\IEEEsetlabelwidth{$P_{i}^{g}/Q_{i}^{g}$}\IEEEusemathlabelsep]
\item[$P_{i}^{g}/Q_{i}^{g}$] Active/reactive power ($\SI{}{\mega\watt}/\SI{}{\mega\VAr}$) generation of generator $g$ at bus $i$.
\item[$P_{ij}/Q_{ij}$] Active/reactive power ($\SI{}{\mega\watt}/\SI{}{\mega\VAr}$) flow along transmission line $ij$.
\item[$V_{i}$] Complex phasor voltage (p.u.) at bus $i$ ($V_{i}=\left|V_{i}\right| \angle \theta_{i}$).
\item[$\theta_{i}$] Voltage angle ($\SI{}{\rad}$) at bus $i$.
\end{IEEEdescription}

\section{Introduction}
The alternating current (AC) power flow equations, which model the steady-state physics of power flows, are the linchpins of a broad spectrum of optimization problems in electrical power systems.
Unfortunately, these nonlinear equations are the main source of nonconvexity in these problems and are notorious for being extremely challenging to solve using global nonlinear programming (GNLP) solvers. Therefore, the research community has focused on improving interior-point based nonlinear optimization methods to compute feasible solutions efficiently \cite{IPMforOPF,MATPOWER}. However, these methods only guarantee local optimality and therefore provide no bounds on the optimal solution.

Due to these challenges, the electricity industry resorted to two main approaches for finding a good tradeoff between computational complexity and quality of lower bound. The first approach consists of methods for approximating the power flow equations, such as the direct current optimal power flow (DC OPF). The DC OPF exploits some physical properties of power flows in typical power systems, such as small bus voltage magnitude ranges and small bus voltage angle differences, to approximate the AC OPF by a linear program (LP). Under normal operating conditions and some adjustments of the lines losses, the DC OPF can approximate the AC active power flow equations with reasonable accuracy \cite{DCOPFrevisited}. Moreover, the DC OPF can be extended to a mixed-integer linear programming (MILP) model to suit a wide variety of optimization applications in power system operations such as optimal transmission switching (OTS), capacitor placement, transmission and distribution network expansion planning, optimal feeder reconfiguration, power system restoration, and vulnerability analysis, to name a few. In summary, the DC OPF is particularly attractive because it leverages the high computational efficiency of LP and MILP solvers. On the downside, the DC OPF fails to capture important physical aspects of power flow, such as the reactive power and voltage magnitude, which are crucial in many applications to ensure voltage stability and AC power flow feasibility. Additionally, the accuracy and feasibility of the DC OPF under congested or unstable operating conditions are questionable. For these reasons, the DC OPF can return solutions that are infeasible in the original space and is proven to be inadequate in applications such as optimal transmission switching \cite{primalanddualboundsforOTS,OTSwithVsecurityandcontingency}.

The second, more recent approach consists of developing computationally efficient convex relaxations. In particular, the second-order cone programming (SOCP) and the semidefinite programming (SDP) relaxations have garnered considerable attention in the electricity industry. The increased interest in this line of research stems from the fact that the SDP relaxation is proven to be exact (i.e. yields a zero optimality gap) on a variety of case studies \cite{0dualitygapinOPF}. However, in many practical OPF instances, the SDP relaxation yields inexact solutions \cite{investigationofnon0gap,InexactnessofSDP}. In these scenarios, an AC feasible solution cannot be recovered from the SDP relaxed solution. The SDP relaxation can be strengthened by solving a hierarchy of moment relaxations at the cost of larger SDP problems \cite{moment-sosopf,Moment-BasedRelaxations}. The main drawback of the SDP relaxation is that it cannot be readily embedded in mixed-integer programming (MIP) models as easily as LP models. Furthermore, mixed-integer SDP technology is still in its infancy compared to the more mature MILP technology. 

Even more recently, increased attention was given to the computationally less demanding SOCP relaxation initially proposed in \cite{radialDNusingConicP}. The SOCP relaxation in its classical form \cite{radialDNusingConicP} is shown to be dominated by the SDP relaxation but recent strengthening techniques \cite{QCrelaxation,strongSOCP,Strengtheningwithboundtightening} have shifted this paradigm. The attractiveness of the SOCP relaxation is also due to the fact that SOCP models can be easily extended to mixed-integer quadratically constrained programming (MIQCP) models to suit applications with discrete variables, mentioned earlier.

Against this background, this paper aims at narrowing the gap between LP approximations and convex relaxations of AC power flow equations by retaining the good lower bounds of the SOCP relaxation and relishing the computational efficiency of LP solvers. In more detail, this paper proposes two LP approximations for the OPF problem based on tight polyhedral approximations of the second-order cone (SOC) constraints \cite{newSOCapproximation}. The first LP model is a direct LP approximation of the classical SOCP relaxation in \cite{radialDNusingConicP}, whereas the second LP model employs strengthening techniques inspired by \cite{QCrelaxation} which aim at preserving stronger links between the voltage variables through convex envelopes of the polar representation. As shown in \cite{QCrelaxation}, a model adopting these strengthening techniques neither dominates nor is dominated by the SDP relaxation. It is important to note that in this context the term ``tight'' designates the high accuracy of the LP approximation of the OPF compared to its respective parent SOCP relaxation.

This paper is not the first attempt to approximate both active and reactive power flow equations in the OPF problem. The LP approximation in \cite{LPOPF} is based on outer approximations which
are strengthened by incorporating several different types of valid inequalities. However, both the computation time and the accuracy of the approximation seem to vary arbitrarily with system size. In contrast to \cite{LPOPF}, the accuracy of the LP models in this paper does not exceed $10^{-2}\%$ in the worst case and the computational efficiency is comparable to, if not better than, that of the SOCP models in most test instances.

In summary, the contributions of this paper are twofold:
\begin{itemize}
\item The two LP models are tested on instances from MATPOWER \cite{MATPOWER} and NESTA v0.5.0 archive \cite{NESTA} with up to $9241$ buses and different operating conditions and are shown to consistently produce high approximation accuracies in the order of $10^{-4}\%$ on average.
\item Numerical results show that the computational efficiency of the LP models is comparable to, if not better than, that of the SOCP models in most instances. This performance is ideal for MILP extensions of these LP models since MILP is computationally more efficient than MIQCP.
\end{itemize}

The paper progresses with the OPF problem formulation in Section~\ref{sec:OPFproblem}, followed by a review of the different types of relaxations proposed in the literature in Section~\ref{sec:Relaxations}. Sections~\ref{sec:Polyhedralformulations} and~\ref{sec:LPOPF} describe the polyhedral formulations of the OPF problem and Section~\ref{sec:evaluation} showcases the numerical results. Finally, the paper concludes in Section~\ref{sec:conclusion}.

\section{Optimal power flow problem}\label{sec:OPFproblem}
In a power network, the OPF problem consists of finding the most economic dispatch of power from generators to satisfy the load at all buses in a way that is governed by physical laws, such as Ohm's Law and Kirchhoff's Law, and other technical restrictions, such as transmission line thermal limit constraints. More specifically, the OPF problem is written as in Model~\ref{Model0}, where $T_{ij}=\tau_{ij}\mathrm{e}^{\mathrm{i} \theta_{ij}^{\text{shift}}}$ is the complex tap ratio of a phase shifting transformer.
\setlength{\belowdisplayskip}{3pt} \setlength{\belowdisplayshortskip}{0pt}
\setlength{\abovedisplayskip}{3pt} \setlength{\abovedisplayshortskip}{0pt}
\small
\begin{Model}
  \DontPrintSemicolon
	\small
\begin{subequations}\label{eq:opf}
		\begin{align}
	 \mbox{\normalsize minimize} & \sum_{(g,i)\in \mathcal{G}} c2^{g}_{i}\left(P_{i}^{g}\right)^2 +c1^{g}_{i}\left(P_{i}^{g}\right)+c0^{g}_{i} & \label{eq:objective}\\
	 \text{\normalsize subject } & \text{\normalsize to}  &\nonumber \\
	\underline{P}_{i}^{g} \leq P_{i}^{g} & \leq \overline{P}_{i}^{g},\ \underline{Q}_{i}^{g} \leq   Q_{i}^{g} \leq \overline{Q}_{i}^{g}, \qquad \ (g,i)\in \mathcal{G} & \label{eq:PQminmax} \\
	\underline{V}_{i} & \leq \left|V_{i}\right| \leq \overline{V}_{i}, \ \ \qquad \qquad \qquad \quad \ \  \ i\in \mathcal{B} & \label{eq:Vminmax} \\ 
	\underline{\theta}_{ij}^{\Delta} & \leq \theta_{i}-\theta_{j} \leq \overline{\theta}_{ij}^{\Delta}, \qquad \qquad \quad \quad \ \  ij \in \mathcal{L} & \label{eq:anglediff} \\
	\sum_{(g,i)\in \mathcal{G}} P_{i}^{g}&-P_{i}^{d} =  \sum_{j\in \mathcal{B}_{i}} P_{ij},\sum_{(g,i)\in \mathcal{G}}  Q_{i}^{g}-Q_{i}^{d}=  \sum_{j\in \mathcal{B}_{i}} Q_{ij}, \nonumber \\ & \qquad \qquad \qquad \qquad \qquad \qquad \qquad \ i\in \mathcal{B}  & \label{eq:KCL}  \\	
	P_{ij}=\Re&\left\{\frac{y_{ij}^*-\mathrm{i}\frac{b^\text{sh}_{ij}}{2}}{\left|T_{ij}\right|^2}\right\} \left|V_{i}\right|^2-\Re\left\{\frac{y_{ij}^*}{T_{ij}}\right\} \Re\left\{V_{i} V_{j}^*\right\} & \nonumber\\
	+\Im&\left\{\frac{y_{ij}^*}{T_{ij}}\right\}\Im\left\{V_{i} V_{j}^*\right\}, \ \ \qquad \qquad \  \ ij \in \mathcal{L} & \label{eq:LF_Pij} \\
	Q_{ij}=\Im&\left\{\frac{y_{ij}^*-\mathrm{i}\frac{b^\text{sh}_{ij}}{2}}{\left|T_{ij}\right|^2}\right\} \left|V_{i}\right|^2-\Im\left\{\frac{y_{ij}^*}{T_{ij}}\right\} \Re\left\{V_{i} V_{j}^*\right\} & \nonumber\\
	-\Re&\left\{\frac{y_{ij}^*}{T_{ij}}\right\}\Im\left\{V_{i} V_{j}^*\right\}, \ \ \qquad \qquad \  \ ij \in \mathcal{L} & \label{eq:LF_Qij} \\
	P_{ji}=\Re&\left\{y_{ji}^*-\mathrm{i}\frac{b^\text{sh}_{ji}}{2}\right\} \left|V_{j}\right|^2-\Re\left\{\frac{y_{ji}^*}{T_{ji}^*}\right\} \Re\left\{V_{j} V_{i}^*\right\} & \nonumber\\
	+\Im&\left\{\frac{y_{ji}^*}{T_{ji}^*}\right\}\Im\left\{V_{j} V_{i}^*\right\}, \ \ \qquad \qquad \  \ ij \in \mathcal{L} & \label{eq:LF_Pji} \\
	Q_{ji}=\Im&\left\{y_{ji}^*-\mathrm{i}\frac{b^\text{sh}_{ji}}{2}\right\} \left|V_{j}\right|^2-\Im\left\{\frac{y_{ji}^*}{T_{ji}^*}\right\} \Re\left\{V_{j} V_{i}^*\right\} & \nonumber\\
	-\Re&\left\{\frac{y_{ji}^*}{T_{ji}^*}\right\}\Im\left\{V_{j} V_{i}^*\right\}, \ \ \qquad \qquad \  \ ij \in \mathcal{L} &\label{eq:LF_Qji} \\
	& \sqrt{P_{ij}^2+Q_{ij}^2} \leq \overline{S}_{ij}, \quad \ \ \ \ \ \ ij \in \mathcal{L} \cup \mathcal{L}_{t}. & \label{eq:linethermallimit} 	
		\end{align}
\end{subequations}		
  \caption{AC OPF}
	\label{Model0}
\end{Model}
\normalsize


Problem~\eqref{eq:opf} is a nonconvex nonlinear optimization problem that is proven to be NP-hard \cite{0dualitygapinOPF}. Therefore, solving large-scale instances of this problem to optimality is intractable. Consequently, applying interior-point methods (IPM) \cite{MATPOWER} to this problem provides no bounds or guarantees on the optimality of the solution, which incited researchers to channel considerable effort on convex relaxation methods.

The next section presents two of the most extensively studied relaxations of problem~\eqref{eq:opf}, namely, the SDP and the SOCP relaxations.

\section{The SDP and SOCP Relaxations}\label{sec:Relaxations}
The SDP relaxation was first introduced in \cite{BaiSDP} and later formalized in \cite{0dualitygapinOPF}. An equivalent formulation of problem \eqref{eq:opf}, described in \cite{0dualitygapinOPF}, starts by setting
\begin{align}
	W = 
 \begin{bmatrix}
  \left|V_{1}\right|^2 & V_{1}V_{2}^* & \cdots & V_{1}V_{\left|\mathcal{B}\right|}^* \\
   V_{2}V_{1}^* & \left|V_{2}\right|^2 & \cdots & V_{2}V_{\left|\mathcal{B}\right|}^* \\
  \vdots  & \vdots  & \ddots & \vdots  \\
  V_{\left|\mathcal{B}\right|}V_{1}^* & V_{\left|\mathcal{B}\right|}V_{2}^* & \cdots & \left|V_{\left|\mathcal{B}\right|}\right|^2 
 \end{bmatrix}
\end{align}
and requiring that $W \succeq 0$ and $\text{rank}(W)=1$. The SDP relaxation is then obtained by dropping the rank constraint. The main setback of applying the SDP relaxation to very large systems is that the matrix $W$ is dense even when all the data matrices are sparse. To this end, sparsity exploiting methods have been proposed in \cite{graphpartitioningtechnique,ExploitingSparsityinSDP,Moment-BasedRelaxations,convexrelaxationmesh} to reduce the computational burden. However, even after applying sparsity exploiting techniques, the computational efficiency of current primal-dual interior-point methods for large-scale SDP is still substantially lower than that of state-of-the-art SOCP solvers. Therefore, in an effort to exploit the sparsity of the power network and leverage the higher computational efficiency of SOCP solvers, \cite{physics} proposes further relaxing some selected positive semidefinite (PSD) conditions in the PSD constraint matrix $W$ to SOC constraints \cite{SOCPrelaxationofSDP}. The first condition is that every $2 \times 2$ principal submatrix of a PSD matrix is also a PSD matrix. The second condition is that the positive semidefiniteness of each $2 \times 2$ symmetric matrix can be represented by a SOC constraint. More specifically, $W \succeq 0$ is replaced by $\left|\mathcal{L}\right|$ constraints of the form
\begin{align}
	\left|W_{ij}\right|^2 \leq W_{ii} W_{jj}, \ \left(W_{ii}, W_{jj} \geq 0\right),  \ \ ij \in \mathcal{L}. \label{eq:SOCconstraint1}
\end{align}
It was also observed in \cite{physics} that the resulting SOCP relaxation \eqref{eq:SOCconstraint1} is tantamount to the SOCP relaxation proposed earlier in \cite{radialDNusingConicP} for radial networks. The SOC representation of the power flow constraints \cref{eq:LF_Pij,eq:LF_Qij,eq:LF_Pji,eq:LF_Qji} was initially introduced in \cite{Reliableloadflow} as follows:
\begin{subequations}
\begin{align}
	W_{ij} & =V_{i} V_{j}^* \\
	W_{ij} W_{ij}^* & =V_{i} V_{j}^* V_{i}^* V_{j}\\
	\left|W_{ij}\right|^2  & = W_{ii} W_{jj}. \label{eq:SOCnonconvex}
\end{align}
\end{subequations}
However, \eqref{eq:SOCnonconvex} is not convex because it describes the surface of a rotated SOC. Therefore a convex relaxation of \eqref{eq:SOCnonconvex} was proposed in \cite{radialDNusingConicP} by relaxing the equality into an inequality as in \eqref{eq:SOCconstraint1}. By defining
\begin{subequations}
\begin{align}
	W_{ii} & =\left|V_{i}\right|^2, \label{eq:Wii} \\
	W_{ij}^\mathrm{r} & =\Re\left\{W_{ij}\right\}=\left|V_{i}\right| \left|V_{j}\right| \cos\left(\theta_{i}-\theta_{j}\right), \label{eq:Wijreal} \\
	W_{ij}^\mathrm{i} & =\Im\left\{W_{ij}\right\}=\left|V_{i}\right| \left|V_{j}\right| \sin\left(\theta_{i}-\theta_{j}\right), \label{eq:Wijimaginary}
\end{align}
\end{subequations}
the SOCP relaxation of problem \eqref{eq:opf} can be written as in Model~\ref{Model1}, where \eqref{eq:anglediff1}, introduced in \cite{convexrelaxationmesh}, is the equivalent of \eqref{eq:anglediff}.
\small
\begin{Model}
  \DontPrintSemicolon
	\small
\begin{subequations}\label{eq:SOCP-0}
		\begin{align}
			\mbox{\normalsize minimize} & \ \sum_{(g,i)\in \mathcal{G}} c2^{g}_{i}\left(P_{i}^{g}\right)^2+c1^{g}_{i}\left(P_{i}^{g}\right)+c0^{g}_{i} & \\
			\text{\normalsize subject }& \text{\normalsize to \small \cref{eq:PQminmax,eq:KCL,eq:linethermallimit,eq:SOCconstraint1},} & \\
			\underline{V}_{i}^2 &\leq W_{ii} \leq \overline{V}_{i}^2, \qquad \qquad \qquad \quad  \ i\in \mathcal{B} & \label{eq:Wminmax} \\
			\tan \left(\underline{\theta}_{ij}^{\Delta}\right)W_{ij}^\mathrm{r} & \leq W_{ij}^\mathrm{i} \leq \tan \left(\overline{\theta}_{ij}^{\Delta}\right) W_{ij}^\mathrm{r}, \qquad   ij \in \mathcal{L} & \label{eq:anglediff1} \\
			P_{ij}=\Re&\left\{\frac{y_{ij}^*-\mathrm{i}\frac{b^\text{sh}_{ij}}{2}}{\left|T_{ij}\right|^2}\right\} W_{ii}-\Re\left\{\frac{y_{ij}^*}{T_{ij}}\right\} W_{ij}^\mathrm{r} & \nonumber\\
	+\Im&\left\{\frac{y_{ij}^*}{T_{ij}}\right\} W_{ij}^\mathrm{i}, \  \qquad \qquad \qquad \  \ ij \in \mathcal{L} & \label{eq:LF_Pij1} \\
	Q_{ij}=\Im&\left\{\frac{y_{ij}^*-\mathrm{i}\frac{b^\text{sh}_{ij}}{2}}{\left|T_{ij}\right|^2}\right\} W_{ii}-\Im\left\{\frac{y_{ij}^*}{T_{ij}}\right\} W_{ij}^\mathrm{r} & \nonumber\\
	-\Re&\left\{\frac{y_{ij}^*}{T_{ij}}\right\}W_{ij}^\mathrm{i}, \  \qquad \qquad \qquad \  \ ij \in \mathcal{L} & \label{eq:LF_Qij1} \\
	P_{ji}=\Re&\left\{y_{ji}^*-\mathrm{i}\frac{b^\text{sh}_{ji}}{2}\right\} W_{jj}-\Re\left\{\frac{y_{ji}^*}{T_{ji}^*}\right\} W_{ij}^\mathrm{r} & \nonumber\\
	-\Im&\left\{\frac{y_{ji}^*}{T_{ji}^*}\right\}W_{ij}^\mathrm{i}, \  \qquad \qquad \qquad \  \ ij \in \mathcal{L} & \label{eq:LF_Pji1} \\
	Q_{ji}=\Im&\left\{y_{ji}^*-\mathrm{i}\frac{b^\text{sh}_{ji}}{2}\right\} W_{jj}-\Im\left\{\frac{y_{ji}^*}{T_{ji}^*}\right\} W_{ij}^\mathrm{r} & \nonumber\\
	+\Re&\left\{\frac{y_{ji}^*}{T_{ji}^*}\right\}W_{ij}^\mathrm{i}, \  \qquad \qquad \qquad \  \ ij \in \mathcal{L}. &\label{eq:LF_Qji1} 	
		\end{align}
\end{subequations}		
  \caption{SOCP-$0$}
	\label{Model1}
\end{Model}
\normalsize

Next, the SOCP relaxation in Model~\ref{Model1} can be strengthened by adding constraints that define tight convex envelopes of the nonlinear terms in \eqref{eq:Wii}, \eqref{eq:Wijreal} and \eqref{eq:Wijimaginary} \cite{convexQR,QCrelaxation}. As shown in \cite{McCormick}, the convex hull of a bilinear term $\left\{w=xy | \left(x,y\right) \in \left[\underline{x},\overline{x}\right] \times\left[\underline{y},\overline{y}\right]\right\}$ is given by
\begin{subnumcases}{\label{eq:w} \hspace{-9mm} \text{conv}\mathcal{M}:= \!}
	\! w \geq \underline{x} y +\underline{y} x -\underline{x}\underline{y} \label{eq:w1} \\
	\! w \geq \overline{x} y +\overline{y} x -\overline{x}\overline{y} \label{eq:w2} \\
	\! w \leq \underline{x} y +\overline{y} x -\underline{x}\overline{y} \label{eq:w3} \\
	\! w \leq \overline{x} y +\underline{y} x -\overline{x}\underline{y} \label{eq:w4},
\end{subnumcases}
and the convex hull of $\left\{w_{2}=x^2 | x \in \left[\underline{x},\overline{x}\right]\right\}$ is given by
\begin{subnumcases}{\label{eq:V2} \hspace{-9mm} \text{conv}\mathcal{C}:= \!}
	\! w_{2} \geq x^2 \label{eq:V21} \\
	\! w_{2} \leq \left(\overline{x}+\underline{x}\right)x-\overline{x} \underline{x} \label{eq:V22}.
\end{subnumcases}

Under the assumption that $\theta^{\Delta}$ does not exceed the range $\left(-\frac{\pi}{2},\frac{\pi}{2}\right)$,\footnote{\label{note1} In practice, $\theta^{\Delta}$ typically does not exceed $\pm 10^{\circ}$ \cite{UsefulnessofDCPF}.} convex envelopes of $\left\{x_{c}=\cos(x) | x \in \left[\underline{x},\overline{x}\right]\right\}$ and $\left\{x_{s}=\sin(x) | x \in \left[\underline{x},\overline{x}\right]\right\}$ are given by
\begin{subnumcases}{\label{eq:Ccosine} \hspace{-9mm} \text{conv}\mathcal{C}_c:= \!}
	\! x_{c} \leq 1-\frac{1-\cos\left(\overline{x}\right)}{\overline{x}^2} x^2 \label{eq:Ccosine1}\\
	\! x_{c} \geq \cos\left(\overline{x}\right). \label{eq:Ccosine2}
\end{subnumcases}

\begin{subnumcases}{\label{eq:Csine} \hspace{-8mm} \text{conv}\mathcal{C}_s:= \!}
	\! x_{s} \leq \cos\left(\frac{\overline{x}}{2}\right) \left(x-\frac{\overline{x}}{2}\right)+\sin\left(\frac{\overline{x}}{2}\right) \label{eq:Csine1} \\
	\! x_{s} \geq \cos\left(\frac{\overline{x}}{2}\right) \left(x+\frac{\overline{x}}{2}\right)-\sin\left(\frac{\overline{x}}{2}\right). \label{eq:Csine2}
\end{subnumcases}

The convex envelopes \cref{eq:w1,eq:w2,eq:w3,eq:w4,eq:V21,eq:V22,eq:Ccosine1,eq:Ccosine2,eq:Csine1,eq:Csine2} are introduced in \cite{convexQR,QCrelaxation} to preserve stronger links between the complex phasor voltage variables. Using these, the strengthened SOCP relaxation of problem~\eqref{eq:opf} is shown in Model~\ref{Model2}.

The next section describes how to tightly approximate Models~\ref{Model1} and~\ref{Model2} by LPs.
\small
\begin{Model}
  \DontPrintSemicolon
	\small
\begin{subequations}\label{eq:SOCP-S}
		\begin{align}
			\mbox{\normalsize minimize} & \ \sum_{(g,i)\in \mathcal{G}} c2^{g}_{i}\left(P_{i}^{g}\right)^2+c1^{g}_{i}\left(P_{i}^{g}\right)+c0^{g}_{i} & \\
			\text{\normalsize subject }& \text{\normalsize to \small \cref{eq:PQminmax,eq:Vminmax,eq:anglediff,eq:KCL}, \cref{eq:linethermallimit,eq:SOCconstraint1}, \cref{eq:Wminmax,eq:anglediff1,eq:LF_Pij1,eq:LF_Qij1,eq:LF_Pji1,eq:LF_Qji1},} & \\
			\text{conv} & \mathcal{C}\left(W_{ii}=\left|V_{i}\right|^2\right), \qquad \qquad \qquad \quad i \in \mathcal{B}& \label{eq:WiiC} \\
			\text{conv} & \mathcal{C}_{c}\left(x_{c,ij}=\cos\left(\theta_{i}-\theta_{j}\right)\right), \qquad \ \ ij \in \mathcal{L}& \label{eq:ccos} \\ 
			\text{conv} & \mathcal{C}_{s}\left(x_{s,ij}=\sin\left(\theta_{i}-\theta_{j}\right)\right), \qquad \ \ ij \in \mathcal{L}& \label{eq:csin} \\
			\text{conv} & \mathcal{M}(w_{ij}=\left|V_{i}\right| \left|V_{j}\right|), \qquad \qquad \quad \ ij \in \mathcal{L}& \label{eq:cVV} \\ 
			\text{conv} & \mathcal{M}(W_{ij}^{\mathrm{r}}=w_{ij} x_{c,ij}), \qquad \qquad \quad ij \in \mathcal{L}& \label{eq:cWcos} \\ 
			\text{conv} & \mathcal{M}(W_{ij}^{\mathrm{i}}=w_{ij} x_{s,ij}), \qquad \qquad \quad ij \in \mathcal{L}.& \label{eq:cWsin}
		\end{align}
\end{subequations}		
  \caption{SOCP-S}
	\label{Model2}
\end{Model}
\normalsize

\section{Polyhedral formulations}\label{sec:Polyhedralformulations}
This section describes how to tightly approximate a 3-dimensional SOC by a polyhedral set. This formulation is also extended to approximate a 4-dimensional rotated SOC. Additionally, this section also describes a polyhedral approximation of the cosine term.

\subsection{Approximation of the 3-dimensional SOC}\label{sec:3DSOC}
A 3-dimensional SOC $\SOC^{2}$ is a subset of $\reals^{3}$ defined by $\SOC^{2}=\left\{\left(r,x_0,y_0\right) \in \reals^{3} |  \sqrt{x_0^2+y_0^2} \leq r\right\}$. One ostensible way to approximate $\SOC^{2}$ is by a regular circumscribed $m$-polyhedral cone $\mathcal{P}_{m} \subseteq \reals^{3}$, which is described by $m$ inequalities. The polyhedron $\mathcal{P}_{m}$ therefore contains $\SOC^{2}$, that is:
\begin{align*}
	\SOC^{2}\subseteq \mathcal{P}_{m} \subseteq \SOC^{2}_{\epsilon}=\left\{\left(r,x_0,y_0\right) \in \reals^{3} |  \sqrt{x_0^2+y_0^2} \leq \left(1+\epsilon \right)r\right\},
\end{align*}
where $\SOC^{2}_{\epsilon}$ is an $\epsilon$-relaxed $\SOC^{2}$ and $\epsilon=\cos(\frac{\pi}{m})^{-1}-1$ is the approximation accuracy. However, this approach requires $233$ linear inequalities even for a relatively modest accuracy of $10^{-4}$. 

As proposed by \cite{PolyhedralApproximationsofSOC}, the key to decreasing the number of inequalities is to lift the approximating polyhedron into a higher dimensional space by introducing several additional variables and projecting it onto the 3-dimensional subspace of the original variables $\left(r,x_0,y_0\right)$. This polyhedral formulation is modified in \cite{newSOCapproximation} to require fewer variables and linear inequalities. In more detail, for an integer $k \geq 2$, let $\mathcal{P}_{k} \in \reals^{2k+3}$ be defined as
\setlength{\belowdisplayskip}{0pt} \setlength{\belowdisplayshortskip}{0pt}
\setlength{\abovedisplayskip}{0pt} \setlength{\abovedisplayshortskip}{0pt}
\[
$$
\begin{subnumcases}{\label{eq:polyhedron} \hspace{-9mm} \mathcal{P}_{k}:= \!}
	\! x_{i+1}=x_{i}\cos(\frac{\pi}{2^i})+y_{i}\sin(\frac{\pi}{2^i}) ,0 \leq i < k, \label{eq:set01}\\
	\! y_{i+1} \geq y_{i}\cos(\frac{\pi}{2^i})-x_{i}\sin(\frac{\pi}{2^i}),0 \leq i < k, \label{eq:set02}\\
	\! y_{i+1} \geq -y_{i}\cos(\frac{\pi}{2^i})+x_{i}\sin(\frac{\pi}{2^i}),0 \leq i < k, \label{eq:set03}\\
	\! \hspace{5mm} r=x_{k}\cos(\frac{\pi}{2^i})+y_{k}\sin(\frac{\pi}{2^i}). \label{eq:set04}
\end{subnumcases}
$$  \nonumber
\]
As identified by \cite{newSOCapproximation}, the projection of set $\mathcal{P}_{k}$ on the subspace of $\left(r,x_0,y_0\right)$ is a polyhedral approximation of $\SOC^{2}$ with accuracy $\epsilon=\cos(\frac{\pi}{2^k})^{-1}-1$.

Polyhedron $\mathcal{P}_{k}$ in~\eqref{eq:polyhedron} requires $2k+3$ variables, $2k$ linear inequalities and $k+1$ linear equalities. This formulation can be reduced even further using \eqref{eq:set01} and \eqref{eq:set04} to replace $x_{i}$ in \eqref{eq:set02} for $i=1,\ldots,k$, and $y_{k}$ in \eqref{eq:set03}. The resulting polyhedron now requires $k+2$ variables $\left(r,x_0,y_0,\ldots,y_{k-1}\right)$ and $2k$ linear inequalities. Now, an accuracy of $10^{-4}$ would require $16$ inequalities and $7$ additional variables (or $10$ variables including the $3$ original ones), i.e. $k=8$. 
 
\subsection{Approximation of the 4-dimensional rotated SOC}\label{sec:rotatedSOC}
\setlength{\belowdisplayskip}{4pt} \setlength{\belowdisplayshortskip}{0pt}
\setlength{\abovedisplayskip}{4pt} \setlength{\abovedisplayshortskip}{0pt}
The formulation in Section~\ref{sec:3DSOC} can be further extended to approximate a 4-dimensional rotated SOC $\SOC^3$ \cite{PolyhedralFormulationsofDNEP}, which is a subset of $\reals^{4}$ defined by $\SOC^3=\left\{\left(r_1,r_2,x_0,y_0\right) \in \reals^{4} |  x_0^2+y_0^2 \leq r_1 r_2\right\}$. The rotated SOC $\SOC^3$ can also be expressed as
\begin{eqnarray}
	\displaystyle{r }&\geq&\displaystyle{ \sqrt{x_0^2+y_0^2}}, \label{eq:r}\\
	r' &\geq& \sqrt{\left(x'_0\right)^2+\left(y'_0\right)^2}, \label{eq:r'}\\
	\displaystyle{r'}&=&\displaystyle{\frac{\left(r_1+r_2\right)}{2}, r=\frac{\left(r_1-r_2\right)}{2}, y'_0=r \label{eq:rcoupling}}.
\end{eqnarray}
Now, \eqref{eq:r} and \eqref{eq:r'} can each be approximated by $\mathcal{P}_{k}$ and coupled by \eqref{eq:rcoupling} to construct a polyhedron $\mathcal{P}_{k}^{r}$ with $2k+6$ variables $\left(r_1,r_2,x_0,y_0,\ldots,y_{k-1},r,r',x'_0,y'_0,\ldots,y'_{k-1}\right)$, $4k$ linear inequalities and $3$ linear equalities to approximate $\SOC^3$.

\subsection{Approximation of the square of a variable}\label{sec:x2}
A function of the form $w2 \geq x^2$ can be approximated by a polyhedron $\mathcal{P}_{k}^{r}$ as described in Section~\ref{sec:rotatedSOC}. However, numerical simulations have shown that the increased accuracy of approximating \eqref{eq:V21} by a polyhedron $\mathcal{P}_{k}^{r}$ has very little effect on the accuracy of the overall solution. This stems from a combination of two factors. The first is the small range of the voltage magnitude (e.g. [0.95,1.05] in practice) and the second is that \eqref{eq:WiiC} might be dominated by \eqref{eq:SOCconstraint1}. To this end, a simpler polyhedral approximation of \eqref{eq:V2} is constructed as follows: $l$ points $x_1,\ldots,x_l$ are selected in the interval $\left[\underline{x},\overline{x}\right]$, which allows adding $l+1$ constraints of the form
\setlength{\belowdisplayskip}{0pt} \setlength{\belowdisplayshortskip}{0pt}
\setlength{\abovedisplayskip}{0pt} \setlength{\abovedisplayshortskip}{0pt}
\[$$
\begin{subnumcases}{\label{eq:x2} \hspace{-6mm} \mathcal{P}^{\text{alt}}_{l}:= \!}
	\! w_2 \geq \left(2 x_{h}\right)x-x_{h}^{2}, \quad h=\left\{1,\ldots,l\right\}, \label{eq:x21}\\
	\! w_{2} \leq \left(\overline{x}+\underline{x}\right)x-\overline{x} \underline{x} . \label{eq:x22}
\end{subnumcases}
$$  \nonumber\]
The approximation in \eqref{eq:x2} requires no additional variables, and numerical simulations have shown that it can result in a high overall accuracy for $l=20$.

\subsection{Approximation of the cosine}
The cosine term in \eqref{eq:Wijreal} can be approximated by a convex affine set provided that $\theta^{\Delta}$ does not exceed the range $\left(-\frac{\pi}{2},\frac{\pi}{2}\right)$\footnotemark[1] \cite{linearapproximationofACflows}. One obvious way is to approximate the quadratic term in \eqref{eq:Ccosine1} by a polyhedron $\mathcal{P}^{\text{alt}}_{l}$ as described in Section~\ref{sec:x2}. However, since the benefit of relaxing the cosine into its convex hull becomes more prominent when $\theta^{\Delta}$ is small, a direct approximation of the cosine can still achieve a high accuracy under these conditions (i.e. when the domain of the cosine function is small). In more detail, a direct polyhedral approximation of the cosine term is constructed as follows:\footnote{\label{note2}The subscript $ij$ is dropped for notational simplicity.} $s$ points $\theta^{\Delta}_1,\ldots,\theta^{\Delta}_s$ are selected in the interval $[\underline{\theta}^{\Delta},\overline{\theta}^{\Delta}]$ and each cosine term $\cos(\theta^{\Delta})$ is replaced with a corresponding new variable $x_{c}$, which allows adding $s+1$ constraints of the form
\[$$
\begin{subnumcases}{\label{eq:cosine} \hspace{-6mm} \mathcal{P}^{\text{cos}}_{s}:= \!}
	\! x_{c} \leq -\sin\left(\theta^{\Delta}_a\right)\left(\theta^{\Delta}-\theta^{\Delta}_a\right) \nonumber \\
	\hspace{8mm} +\cos\left(\theta^{\Delta}_a\right), \quad a=\left\{1,\ldots,s\right\}, \label{eq:cos1}\\
	\! x_{c} \geq \cos\left(\overline{\theta}^{\Delta}\right). \label{eq:cos2}
\end{subnumcases}
$$  \nonumber\]
The approximation in \eqref{eq:cosine} requires no additional variables, and numerical simulations have shown that it can result in a high overall accuracy for $s=20$.

\section{LP Optimal Power Flow}\label{sec:LPOPF}
Given the building blocks in Section~\ref{sec:Polyhedralformulations}, a tight LP approximation of the OPF problem is now possible. The only remaining step is to substitute the quadratic terms in the objective function by corresponding variables and rotated SOC constraints. This substitution now enables leveraging the techniques in Section~\ref{sec:Polyhedralformulations} to tightly approximate the quadratic terms in the objective function by polyhedrons and thereby obtaining a LP approximation of the OPF problem. Specifically, $\left|\mathcal{G}\right|$ variables and constraints of the form 
\setlength{\belowdisplayskip}{4pt} \setlength{\belowdisplayshortskip}{0pt}
\setlength{\abovedisplayskip}{4pt} \setlength{\abovedisplayshortskip}{0pt}
\begin{align}
	p_{g}=\sqrt{c2^{g}_{i}}P^{g}_{i}, \qquad \qquad  \quad (g,i) \in \mathcal{G}, \label{eq:pcoupling}
\end{align}
are introduced along with $N=\left\lfloor \left|\mathcal{G}\right| /2\right\rfloor+\left\lceil \left|\mathcal{G}\right| /2 - \left\lfloor \left|\mathcal{G}\right| /2\right\rfloor \right\rceil$ variables $\alpha_{n}$ and constraints of the form
\begin{align}
	\alpha_{n} \geq p_{2n-1}^2+p_{2n}^2, \qquad \ n \in \left\{1,\ldots,N\right\}.
\end{align}
Finally, the LP approximations of SOCP-$0$ and SOCP-S are shown in Models~\ref{Model3} and~\ref{Model4} respectively, and their accuracy and computational efficiency are evaluated in the next section.
\setlength{\belowdisplayskip}{0pt} \setlength{\belowdisplayshortskip}{0pt}
\setlength{\abovedisplayskip}{0pt} \setlength{\abovedisplayshortskip}{0pt}
\small
\begin{Model}[t!]
  \DontPrintSemicolon
	\small
\begin{subequations}\label{eq:LP0}
		\begin{align}
			\mbox{\normalsize minimize} & \ \sum_{n=1}^{N} \alpha_{n} + \sum_{(g,i)\in \mathcal{G}} c1^{g}_{i}\left(P_{i}^{g}\right)+c0^{g}_{i} & \\
			\text{\normalsize subject }& \text{\normalsize to \small \cref{eq:PQminmax,eq:KCL}, \cref{eq:Wminmax,eq:anglediff1,eq:LF_Pij1,eq:LF_Qij1,eq:LF_Pji1,eq:LF_Qji1}, \eqref{eq:pcoupling},} & \label{eq:LP0const} \\
			\mathcal{P}_{k}^r & \left( pg_{2n-1}^2+pg_{2n}^2 \leq \alpha_{n} \right), \ \ n \in \left\{1,\ldots,N\right\} \label{eq:Pobj} \\
			\mathcal{P}_{k}^r & \left((W_{ij}^{\mathrm{r}})^2 +(W_{ij}^{\mathrm{i}})^2 \leq W_{ii} W_{jj}\right), \ \  ij \in \mathcal{L} \label{eq:PW} \\
			\mathcal{P}_{k} & \left(\sqrt{P_{ij}^2+Q_{ij}^2} \leq \overline{S}_{ij}\right), \qquad \ \  ij \in \mathcal{L} \cup \mathcal{L}_{t}. \label{eq:linethermallimitLP}	 
		\end{align}
\end{subequations}		
  \caption{LP-$0$}
	\label{Model3}
\end{Model}
\begin{Model}[t!]
  \DontPrintSemicolon
	\small
\begin{subequations}\label{eq:LPS}
		\begin{align}
			\mbox{\normalsize minimize} & \ \sum_{n=1}^{N} \alpha_{n} + \sum_{(g,i)\in \mathcal{G}} c1^{g}_{i}\left(P_{i}^{g}\right)+c0^{g}_{i} & \\
			\hspace{-8mm} \text{\normalsize subject }\text{\normalsize to } & \text{\small \cref{eq:LP0const}, \cref{eq:Vminmax,eq:anglediff}, \cref{eq:csin,eq:cVV,eq:cWcos,eq:cWsin}, \cref{eq:Pobj,eq:PW,eq:linethermallimitLP},} & \\
			\mathcal{P}_{l}^{\text{alt}} & \left(W_{ii}=\left|V_{i}\right|^2 \right), \qquad \ \ \qquad \qquad \quad i \in \mathcal{B} \\
			\mathcal{P}_{s}^{\text{cos}} & \left(x_{c,ij}=\cos\left(\theta_{i}-\theta_{j}\right)\right), \qquad \quad \ \ ij \in \mathcal{L}.
		\end{align}
\end{subequations}		
  \caption{LP-S}
	\label{Model4}
\end{Model}
\normalsize

\section{Numerical Evaluation}\label{sec:evaluation}
This section evaluates the accuracy and computational efficiency of the LP approximations in Models~\ref{Model3} and~\ref{Model4} as compared to their respective parent SOCP relaxations in Models~\ref{Model1} and~\ref{Model2}. The models are tested on standard IEEE instances available from the IPM-based OPF solver MATPOWER \cite{MATPOWER} as well as more challenging instances from NESTA v0.5.0 archive \cite{NESTA}, with $k=16$ for the LP models.\footnote{For $k=16$, $\mathcal{P}_{k}$ approximates $\SOC^{2}$ with accuracy $\epsilon=1.15 \times 10^{-9}$.} All simulations are carried out on an Intel Core i7, 3.70GHz, 64-bit, 128GB RAM computing platform. MATPOWER is used to solve the original nonconvex AC model in problem~\eqref{eq:opf}, which provides an upper bound on the optimal solution. Additionally, IPOPT \cite{IPOPT} via the MATLAB toolbox OPTI \cite{OPTI} is used to compute upper bounds on instances where MATPOWER diverges or fails to compute a solution. These locally optimal solutions are not shown in this paper due to space limitations. Interested readers are instead referred to \cite{NESTA} for a complete list of AC (locally optimal) solutions with the exception of MATPOWER's case1354pegase and case9241pegase whose solutions are $\SI{}{\$}74069.354$ and $\SI{}{\$}315912.848$ respectively.\footnote{Computed using IPOPT via OPTI.}

On the other hand, both CPLEX 12.6 \cite{cplex} and Gurobi 6.0.5 \cite{gurobi} are considered for solving the SOCP and the LP models. An interesting observation is that the polyhedral approximations described in Sections~\ref{sec:3DSOC} and~\ref{sec:rotatedSOC}, which are the cornerstones of LP-$0$ and LP-S, make these models particularly difficult to solve using the primal or dual simplex methods. This could be due to the large coefficient ranges and/or due to the irrational coefficients ($\cos(\pi / 2^i)$ and $\sin(\pi / 2^i)$) introduced by these polyhedral formulations. For the LP models, both CPLEX and Gurobi use their default concurrent optimization algorithms which invoke multiple methods (primal simplex, dual simplex and parallel barrier) simultaneously on multiple cores, and return the optimal solution from the method that finishes first. Therefore, in this scenario, only the parallel barrier method is chosen instead of the default concurrent optimization algorithm to solve the LP models. Ultimately, CPLEX is chosen to solve the LP models due to a better performance of its parallel barrier method, for these specific LP models, as compared to Gurobi's parallel barrier method.
\begin{table}[!t]
\centering
\renewcommand{\arraystretch}{1.3}
\caption{Model comparison on MATPOWER instances.}
\resizebox{\linewidth}{!}{%
\begin{tabular}{| c | c | c | c | c | c | c | c | c |}
\cline{2-9}
\multicolumn{1}{c|}{} & \multicolumn{4}{c|}{Optimality Gap (\%)} & \multicolumn{4}{c|}{CPU time ($\SI{}{\second}$)} \\\hline
Case   & SOCP-$0$  & LP-0  & SOCP-S  & LP-S  & SOCP-$0$  & LP-0  & SOCP-S  & LP-S     \\\hline
118	&	0.25	&	0.25	&	0.25	&	0.25	&	0.19	&	0.87	&	0.75	&	1.56	\\\hline
300	&	0.15	&	0.15	&	0.15	&	0.15	&	0.50	&	1.70	&	1.78	&	4.99	\\\hline
1354	&	0.09*	&	0.08	&	0.09*	&	0.08	&	7.72	&	5.76	&	15.08	&	11.93	\\\hline
3375wp	&	0.27*	&	0.26	&	0.26*	&	0.25	&	21.61	&	19.24	&	49.61	&	55.41	\\\hline
9241	&	2.02*	&	2.01	&	2.02*	&	2.01	&	67.44	&	142.20	&	335.53	&	284.47	\\\hline
\end{tabular}}
\label{OPFmatpower}
\end{table}
\begin{table}[!t]
\centering
\renewcommand{\arraystretch}{1.3}
\caption{Model comparison on NESTA instances.}
\resizebox{\linewidth}{!}{%
\begin{tabular}{| c | c | c | c | c | c | c | c | c |}
\cline{2-9}
\multicolumn{1}{c|}{} & \multicolumn{4}{c|}{Optimality Gap (\%)} & \multicolumn{4}{c|}{CPU time ($\SI{}{\second}$)} \\\hline
Case   & SOCP-$0$  & LP-0  & SOCP-S  & LP-S  & SOCP-$0$  & LP-0  & SOCP-S  & LP-S     \\\hline
\multicolumn{9}{|c|}{Normal Operating Conditions} \\\hline
24	&	0.01	&	0.01	&	0.01	&	0.01	&	0.11	&	0.19	&	0.22	&	0.28	\\\hline
29	&	0.14	&	0.14	&	0.12	&	0.12	&	0.30	&	0.28	&	1.56	&	0.69	\\\hline
30\_as	&	0.06	&	0.06	&	0.06	&	0.06	&	0.06	&	0.19	&	0.17	&	0.31	\\\hline
30\_fsr	&	0.39	&	0.39	&	0.39	&	0.39	&	0.06	&	0.17	&	0.19	&	0.25	\\\hline
30	&	15.88	&	15.88	&	15.64*	&	15.62	&	0.06	&	0.17	&	0.20	&	0.19	\\\hline
39	&	0.05	&	0.05	&	0.05	&	0.05	&	0.13	&	0.28	&	0.37	&	0.50	\\\hline
57	&	0.07*	&	0.06	&	0.07*	&	0.06	&	0.13	&	0.20	&	0.51	&	0.30	\\\hline
73	&	0.03	&	0.03	&	0.03	&	0.03	&	0.22	&	0.44	&	1.53	&	1.06	\\\hline
89	&	0.17	&	0.17	&	0.17	&	0.17	&	0.76	&	0.90	&	6.55	&	2.76	\\\hline
118	&	2.07	&	2.07	&	1.72	&	1.72	&	0.45	&	0.61	&	1.97	&	0.72	\\\hline
162	&	4.10*	&	4.03	&	4.00	&	4.00	&	0.62	&	0.70	&	1.78	&	1.34	\\\hline
189	&	0.23	&	0.23	&	0.22	&	0.22	&	2.71	&	2.00	&	2.26	&	1.98	\\\hline
300	&	1.19*	&	1.18	&	1.18	&	1.18	&	1.40	&	1.61	&	4.30	&	2.54	\\\hline
1354	&	0.10*	&	0.08	&	0.09*	&	0.08	&	8.58	&	7.19	&	23.88	&	10.34	\\\hline
2383wp	&	1.08*	&	1.05	&	1.06*	&	1.04	&	14.79	&	21.73	&	45.41	&	34.87	\\\hline
2869	&	0.10*	&	0.09	&	0.10*	&	0.09	&	23.23	&	28.77	&	79.50	&	52.21	\\\hline
3012wp	&	1.06*	&	1.02	&	1.04*	&	1.01	&	18.69	&	19.13	&	43.45	&	59.03	\\\hline
3120sp	&	0.58*	&	0.55	&	0.58*	&	0.60*	&	19.19	&	21.28	&	44.29	&	41.42	\\\hline
3375wp	&	0.52	&	0.52	&	0.51	&	0.51	&	21.09	&	23.21	&	58.00	&	41.99	\\\hline
9241	&	1.76*	&	1.75	&	1.68*	&	1.67	&	356.21	&	112.74	&	588.14	&	486.04	\\\hline
\multicolumn{9}{|c|}{Congested Operating Conditions} \\\hline
14	&	1.34	&	1.34	&	1.34	&	1.34	&	0.05	&	0.06	&	0.06	&	0.08	\\\hline
29	&	0.44	&	0.44	&	0.43	&	0.42	&	0.26	&	0.22	&	1.39	&	0.61	\\\hline
30\_as	&	4.76	&	4.76	&	4.76	&	4.76	&	0.08	&	0.14	&	0.22	&	0.20	\\\hline
30\_fsr	&	45.97	&	45.97	&	45.97	&	45.97	&	0.06	&	0.13	&	0.19	&	0.22	\\\hline
30	&	1.01	&	1.01	&	1.01	&	1.01	&	0.08	&	0.13	&	0.16	&	0.22	\\\hline
39	&	2.99	&	2.99	&	2.97	&	2.97	&	0.09	&	0.14	&	0.42	&	0.33	\\\hline
57	&	0.21	&	0.21	&	0.21	&	0.21	&	0.13	&	0.39	&	0.59	&	0.26	\\\hline
73	&	14.34	&	14.34	&	12.01*	&	12.00	&	0.23	&	0.25	&	2.08	&	0.48	\\\hline
89	&	20.44	&	20.43	&	20.39	&	20.39	&	0.78	&	0.92	&	7.61	&	2.00	\\\hline
118	&	44.08	&	44.08	&	43.93	&	43.93	&	0.38	&	0.48	&	1.83	&	0.78	\\\hline
162	&	1.50*	&	1.34	&	1.33	&	1.33	&	0.78	&	0.80	&	1.84	&	1.42	\\\hline
189	&	6.45*	&	5.79	&	5.84*	&	5.79	&	0.70	&	1.42	&	2.15	&	2.47	\\\hline
300	&	0.84	&	0.84	&	0.82	&	0.82	&	1.22	&	1.51	&	2.84	&	1.92	\\\hline
1354	&	0.58*	&	0.56	&	0.56*	&	0.55	&	11.34	&	6.30	&	31.12	&	24.70	\\\hline
2383wp	&	1.12	&	1.12	&	1.12	&	1.12	&	17.85	&	20.72	&	47.36	&	26.04	\\\hline
2869	&	1.50*	&	1.49	&	1.49	&	1.49	&	32.74	&	30.55	&	108.53	&	67.72	\\\hline
3012wp	&	0.90	&	0.90	&	0.89	&	0.89	&	42.76	&	21.75	&	87.69	&	83.77	\\\hline
3120sp	&	3.03	&	3.03	&	3.03*	&	3.01	&	46.88	&	26.88	&	87.67	&	24.27	\\\hline
3375wp	&	0.59	&	0.60	&	0.59	&	0.59	&	46.41	&	21.29	&	135.38	&	53.54	\\\hline
9241	&	2.59	&	2.59	&	2.46*	&	2.44	&	345.82	&	158.42	&	855.98	&	499.27	\\\hline
\multicolumn{9}{|c|}{Small Angle Difference Conditions} \\\hline
14	&	0.06	&	0.06	&	0.06	&	0.06	&	0.05	&	0.06	&	0.06	&	0.09	\\\hline
24	&	11.42	&	11.42	&	3.88	&	3.88	&	0.13	&	0.20	&	0.25	&	0.39	\\\hline
29	&	34.47	&	34.47	&	20.58*	&	20.57	&	0.30	&	0.52	&	1.84	&	1.15	\\\hline
30\_as	&	9.16	&	9.16	&	3.07	&	3.07	&	0.08	&	0.19	&	0.20	&	0.26	\\\hline
30\_fsr	&	0.62	&	0.62	&	0.56	&	0.56	&	0.08	&	0.17	&	0.22	&	0.26	\\\hline
30	&	5.84	&	5.84	&	3.96	&	3.96	&	0.08	&	0.16	&	0.17	&	0.25	\\\hline
39	&	0.11	&	0.11	&	0.04	&	0.04	&	0.13	&	0.26	&	0.36	&	0.47	\\\hline
57	&	0.11	&	0.11	&	0.10	&	0.10	&	0.13	&	0.23	&	0.61	&	0.31	\\\hline
73	&	8.37	&	8.37	&	3.51	&	3.51	&	0.28	&	0.45	&	2.14	&	1.05	\\\hline
89	&	0.29*	&	0.28	&	0.19*	&	0.18	&	0.81	&	0.72	&	6.66	&	1.18	\\\hline
118	&	12.89	&	12.89	&	8.32	&	8.32	&	0.47	&	0.44	&	2.15	&	0.78	\\\hline
162	&	7.12*	&	7.08	&	6.91	&	6.91	&	0.70	&	0.95	&	1.51	&	1.26	\\\hline
189	&	2.25	&	2.25	&	2.32*	&	2.22	&	1.09	&	2.04	&	2.92	&	1.43	\\\hline
300	&	1.26	&	1.26	&	1.16	&	1.16	&	1.64	&	1.12	&	3.79	&	2.48	\\\hline
1354	&	0.10*	&	0.08	&	0.08*	&	0.07	&	9.34	&	5.32	&	31.17	&	11.34	\\\hline
2383wp	&	4.02	&	4.02	&	3.00*	&	2.97	&	19.13	&	36.64	&	65.72	&	29.23	\\\hline
2869	&	0.16*	&	0.15	&	0.15*	&	0.14	&	27.80	&	19.61	&	80.93	&	44.16	\\\hline
3012wp	&	2.16*	&	2.12	&	1.97*	&	1.92	&	21.96	&	22.93	&	76.36	&	39.89	\\\hline
3120sp	&	2.82*	&	2.79	&	2.62*	&	2.57	&	25.19	&	20.73	&	101.29	&	37.11	\\\hline
3375wp	&	0.53*	&	0.52	&	0.49*	&	0.48	&	24.01	&	24.26	&	91.34	&	48.45	\\\hline
9241	&	1.76*	&	1.75	&	0.81*	&	0.80	&	321.02	&	161.91	&	615.88	&	598.28	\\\hline
\end{tabular}}
\label{OPFnesta}
\end{table}

By letting $S^{AC}$ denote the best known AC solution and $S^{\text{conv}}$ denote the solution from the corresponding relaxation, the optimality gap can be measured as $\left(S^{AC}-S^{\text{conv}}/S^{AC}\right) \times 100$. The optimality gaps and the computation times of the four models are summarized in Tables~\ref{OPFmatpower} and~\ref{OPFnesta} for MATPOWER and NESTA instances respectively. It is evident from Tables~\ref{OPFmatpower} and~\ref{OPFnesta} that both LP-$0$ and LP-S tightly approximate their parent SOCP models, SOCP-$0$ and SOCP-S respectively. However, the values marked by * designate instances where the SOCP relaxation's solution does not match the LP one despite the ``optimal'' exitflag or vice versa. In these cases (*), both Gurobi and IPOPT are used to ascertain that the LP solution is in fact the accurate one in most cases. This is also corroborated by results in the literature, namely in \cite{NESTA} for SOCP-$0$ and \cite{QCrelaxation} for SOCP-S. These discrepancies are due to numerical stability issues despite the solver reporting reaching an optimal solution. This is not surprising since it was also pointed out in \cite{QCrelaxation} that IPOPT is numerically more stable than both CPLEX and Gurobi's QCP for large systems, and was ultimately used for solving their SOCP models. However, CPLEX is still used to solve the SOCP models in this paper for the sake of comparison. Also, the fact that Gurobi and CPLEX are both state-of-the-art LP (and MILP) solvers, it would not make sense to use IPOPT to solve the LP models. In fact, the approximation accuracy is in the order $10^{-5} \%$ when the solution of both SOCP models and LP models does not run into numerical stability issues.

Moreover, Tables~\ref{OPFmatpower} and~\ref{OPFnesta} also show that the computational efficiency of the LP models is comparable to, if not better than, that of the SOCP models in most cases. This performance is ideal for MILP extensions of these LP models, which gives more edge over the MIQCP extensions of the SOCP models because state-of-the-art MIQCP technology is still not as mature as state-of-the-art MILP technology.

\section{Conclusion}\label{sec:conclusion}
Two tight LP approximations of the OPF problem, founded on tight polyhedral approximations of the SOC constraints, are proposed in this paper. The first LP model is a direct LP approximation of the classical SOCP relaxation whereas the second LP model employs strengthening techniques that preserve stronger links between the voltage variables through convex envelopes of the polar representation. Rigorous computational tests on systems with up to $9241$ buses and different operating conditions have shown that the proposed LP models consistently produce high approximation accuracies of $10^{-4}\%$ on average compared to their respective parent SOCP relaxations. Moreover, the computational efficiency of the two proposed LP models is shown to be comparable to, if not better than, that of the SOCP models in most instances, which makes them ideal for MILP extensions knowing that MILP technology is more mature than the MIQCP technology. Finally, the LP models in this paper can easily be extended to handle any convex generator cost function. 

\bibliographystyle{IEEEtran}
{\footnotesize
\bibliography{LPforOPF}}

\begin{thebibliography}{10}
\providecommand{\url}[1]{#1}
\csname url@samestyle\endcsname
\providecommand{\newblock}{\relax}
\providecommand{\bibinfo}[2]{#2}
\providecommand{\BIBentrySTDinterwordspacing}{\spaceskip=0pt\relax}
\providecommand{\BIBentryALTinterwordstretchfactor}{4}
\providecommand{\BIBentryALTinterwordspacing}{\spaceskip=\fontdimen2\font plus
\BIBentryALTinterwordstretchfactor\fontdimen3\font minus
  \fontdimen4\font\relax}
\providecommand{\BIBforeignlanguage}[2]{{%
\expandafter\ifx\csname l@#1\endcsname\relax
\typeout{** WARNING: IEEEtran.bst: No hyphenation pattern has been}%
\typeout{** loaded for the language `#1'. Using the pattern for}%
\typeout{** the default language instead.}%
\else
\language=\csname l@#1\endcsname
\fi
#2}}
\providecommand{\BIBdecl}{\relax}
\BIBdecl

\bibitem{IPMforOPF}
R.~Jabr, A.~Coonick, and B.~Cory, ``A primal-dual interior point method for
  optimal power flow dispatching,'' \emph{Power Systems, IEEE Transactions on},
  vol.~17, no.~3, pp. 654--662, Aug 2002.

\bibitem{MATPOWER}
R.~Zimmerman, C.~Murillo-S\'{a}nchez, and R.~Thomas, ``{MATPOWER}: Steady-state
  operations, planning, and analysis tools for power systems research and
  education,'' \emph{Power Systems, IEEE Transactions on}, vol.~26, no.~1, pp.
  12--19, Feb 2011.

\bibitem{DCOPFrevisited}
B.~Stott, J.~Jardim, and O.~Alsac, ``{DC} power flow revisited,'' \emph{Power
  Systems, IEEE Transactions on}, vol.~24, no.~3, pp. 1290--1300, Aug 2009.

\bibitem{primalanddualboundsforOTS}
C.~Coffrin, H.~Hijazi, K.~Lehmann, and P.~Van~Hentenryck, ``Primal and dual
  bounds for optimal transmission switching,'' in \emph{Power Systems
  Computation Conference (PSCC), 2014}, Aug 2014, pp. 1--8.

\bibitem{OTSwithVsecurityandcontingency}
M.~Khanabadi, H.~Ghasemi, and M.~Doostizadeh, ``Optimal transmission switching
  considering voltage security and {N}-1 contingency analysis,'' \emph{Power
  Systems, IEEE Transactions on}, vol.~28, no.~1, pp. 542--550, Feb 2013.

\bibitem{0dualitygapinOPF}
J.~Lavaei and S.~Low, ``Zero duality gap in optimal power flow problem,''
  \emph{Power Systems, IEEE Transactions on}, vol.~27, no.~1, pp. 92--107, Feb
  2012.

\bibitem{investigationofnon0gap}
D.~K. Molzahn, B.~C. Lesieutre, and C.~L. DeMarco, ``Investigation of non-zero
  duality gap solutions to a semidefinite relaxation of the optimal power flow
  problem,'' in \emph{System Sciences (HICSS), 2014 47th Hawaii International
  Conference on}.\hskip 1em plus 0.5em minus 0.4em\relax IEEE, 2014, pp.
  2325--2334.

\bibitem{InexactnessofSDP}
B.~Kocuk, S.~S. Dey, and X.~A. Sun, ``Inexactness of sdp relaxation and valid
  inequalities for optimal power flow,'' \emph{IEEE Transactions on Power
  Systems}, vol.~31, no.~1, pp. 642--651, Jan 2016.

\bibitem{moment-sosopf}
C.~Josz, J.~Maeght, P.~Panciatici, and J.~Gilbert, ``Application of the
  moment-{SOS} approach to global optimization of the {OPF} problem,''
  \emph{Power Systems, IEEE Transactions on}, vol.~30, no.~1, pp. 463--470, Jan
  2015.

\bibitem{Moment-BasedRelaxations}
D.~Molzahn and I.~Hiskens, ``Sparsity-exploiting moment-based relaxations of
  the optimal power flow problem,'' \emph{Power Systems, IEEE Transactions on},
  vol.~30, no.~6, pp. 3168--3180, Nov 2015.

\bibitem{radialDNusingConicP}
R.~Jabr, ``Radial distribution load flow using conic programming,'' \emph{Power
  Systems, IEEE Transactions on}, vol.~21, no.~3, pp. 1458--1459, Aug 2006.

\bibitem{QCrelaxation}
C.~Coffrin, H.~L. Hijazi, and P.~Van~Hentenryck, ``The {QC} relaxation:
  Theoretical and computational results on optimal power flow,'' \emph{arXiv
  preprint arXiv:1502.07847}, 2015.

\bibitem{strongSOCP}
B.~Kocuk, S.~S. Dey, and X.~A. Sun, ``Strong {SOCP} relaxations for the optimal
  power flow problem,'' \emph{Operations Research}, 2016.

\bibitem{Strengtheningwithboundtightening}
C.~Coffrin, H.~Hijazi, and P.~Van~Hentenryck, ``Strengthening convex
  relaxations with bound tightening for power network optimization,'' in
  \emph{Principles and Practice of Constraint Programming}, ser. Lecture Notes
  in Computer Science, G.~Pesant, Ed.\hskip 1em plus 0.5em minus 0.4em\relax
  Springer International Publishing, 2015, vol. 9255, pp. 39--57.

\bibitem{newSOCapproximation}
F.~Glineur, ``{Topics in Convex Optimization: Interior-Point Methods, Conic
  Duality and Approximations},'' Theses, {Polytechnic College of Mons}, Jan.
  2001.

\bibitem{LPOPF}
D.~Bienstock and G.~Munoz, ``On linear relaxations of {OPF} problems,''
  \emph{arXiv preprint arXiv:1411.1120}, 2014.

\bibitem{NESTA}
\BIBentryALTinterwordspacing
C.~Coffrin, D.~Gordon, and P.~Scott, ``{NESTA}, the {NICTA} energy system test
  case archive,'' \emph{CoRR}, vol. abs/1411.0359, 2014. [Online]. Available:
  \url{http://arxiv.org/abs/1411.0359.}
\BIBentrySTDinterwordspacing

\bibitem{BaiSDP}
X.~Bai, H.~Wei, K.~Fujisawa, and Y.~Wang, ``Semidefinite programming for
  optimal power flow problems,'' \emph{International Journal of Electrical
  Power \& Energy Systems}, vol.~30, no. 6â€“7, pp. 383 -- 392, 2008.

\bibitem{graphpartitioningtechnique}
X.~Bai and H.~Wei, ``A semidefinite programming method with graph partitioning
  technique for optimal power flow problems,'' \emph{International Journal of
  Electrical Power \& Energy Systems}, vol.~33, no.~7, pp. 1309 -- 1314, 2011.

\bibitem{ExploitingSparsityinSDP}
R.~Jabr, ``Exploiting sparsity in {SDP} relaxations of the {OPF} problem,''
  \emph{Power Systems, IEEE Transactions on}, vol.~27, no.~2, pp. 1138--1139,
  May 2012.

\bibitem{convexrelaxationmesh}
R.~Madani, S.~Sojoudi, and J.~Lavaei, ``Convex relaxation for optimal power
  flow problem: Mesh networks,'' \emph{Power Systems, IEEE Transactions on},
  vol.~30, no.~1, pp. 199--211, Jan 2015.

\bibitem{physics}
S.~Sojoudi and J.~Lavaei, ``Physics of power networks makes hard optimization
  problems easy to solve,'' in \emph{Power and Energy Society General Meeting,
  2012 IEEE}, July 2012, pp. 1--8.

\bibitem{SOCPrelaxationofSDP}
S.~Kim, M.~Kojima, and M.~Yamashita, ``Second order cone programming relaxation
  of a positive semidefinite constraint,'' \emph{Optimization Methods and
  Software}, vol.~18, no.~5, pp. 535--541, 2003.

\bibitem{Reliableloadflow}
A.~Gomez~Esposito and E.~Ramos, ``Reliable load flow technique for radial
  distribution networks,'' \emph{Power Systems, IEEE Transactions on}, vol.~14,
  no.~3, pp. 1063--1069, Aug 1999.

\bibitem{convexQR}
H.~Hijazi, C.~Coffrin, and P.~Van~Hentenryck, ``Convex quadratic relaxations of
  mixed-integer nonlinear programs in power systems,'' \emph{Published online
  at http://www. optimization-online. org/DB\_ HTML/2013/09/4057. html}, 2013.

\bibitem{McCormick}
G.~McCormick, ``Computability of global solutions to factorable nonconvex
  programs: Part {I} - convex underestimating problems,'' \emph{Mathematical
  Programming}, vol.~10, no.~1, pp. 147--175, 1976.

\bibitem{UsefulnessofDCPF}
K.~Purchala, L.~Meeus, D.~Van~Dommelen, and R.~Belmans, ``Usefulness of {DC}
  power flow for active power flow analysis,'' in \emph{Power Engineering
  Society General Meeting, 2005. IEEE}, June 2005, pp. 454--459 Vol. 1.

\bibitem{PolyhedralApproximationsofSOC}
A.~Ben-Tal and A.~Nemirovski, ``On polyhedral approximations of the
  second-order cone,'' \emph{Mathematics of Operations Research}, vol.~26,
  no.~2, pp. pp. 193--205, 2001.

\bibitem{PolyhedralFormulationsofDNEP}
R.~A. Jabr, ``Polyhedral formulations and loop elimination constraints for
  distribution network expansion planning,'' \emph{IEEE Transactions on Power
  Systems}, vol.~28, no.~2, pp. 1888--1897, 2013.

\bibitem{linearapproximationofACflows}
C.~Coffrin and P.~Van~Hentenryck, ``A linear-programming approximation of {AC}
  power flows,'' \emph{INFORMS Journal on Computing}, vol.~26, no.~4, pp.
  718--734, 2014.

\bibitem{IPOPT}
A.~W\"{a}chter and L.~T. Biegler, ``On the implementation of an interior-point
  filter line-search algorithm for large-scale nonlinear programming,''
  \emph{Mathematical Programming}, vol. 106, no.~1, pp. 25--57, 2006.

\bibitem{OPTI}
J.~Currie and D.~I. Wilson, ``{OPTI: Lowering the Barrier Between Open Source
  Optimizers and the Industrial MATLAB User},'' in \emph{{Foundations of
  Computer-Aided Process Operations}}, N.~Sahinidis and J.~Pinto, Eds.,
  Savannah, Georgia, USA, 8--11 January 2012.

\bibitem{cplex}
{IBM ILOG CPLEX Optimizer v12.6}.

\bibitem{gurobi}
{Gurobi Optimization Inc.}, ``Gurobi optimizer reference manual,'' 2015.

\end{thebibliography}

\end{document}